\documentclass{article}

\usepackage{amsmath}
\usepackage{amsthm}
\usepackage{epsfig}
\usepackage{graphics}
\usepackage{amsfonts}
\usepackage{amssymb}

\title{Higher Order Birkhoff Averages\footnote{To appear in Dynamical Systems: An International Journal}}

\author{
Thomas Jordan\\
Department of Mathematics\\
University of Bristol\\
Bristol, UK\\
{\small thomas.jordan@bristol.ac.uk}
\and
Vincent Naudot\\
Department of Mathematics\\
Florida Atlantic University\\
Boca Raton, FL, USA\\
{\small  vnaudot@fau.edu }
\and
Todd Young\footnote{Corresponding author.} \\
Department of Mathematics\\
Ohio University\\
Athens, Ohio, USA\\
{\small young@math.ohiou.edu}
}

\numberwithin{equation}{section}

\newcommand{\D}[1]{{\mathbb#1}}
\newcommand{\NN}{{\D{N}}}
\newcommand{\RR}{{\D{R}}}

\def\ep{\epsilon}
\def\la{\lambda}

\newtheorem{defn}{Definition}[section]

\newtheorem{thm}[defn]{Theorem}

\newtheorem{cor}[defn]{Corollary}
\newtheorem{lem}[defn]{Lemma}
\newtheorem{prop}[defn]{Proposition}

\begin{document}

\maketitle

\footnotetext{{\bf AMS Subject Classifcation:} 37A99, 37B10, 37C29}

\begin{abstract}
There are well-known examples of dynamical systems for which the Birkhoff averages with
respect to a given observable along some or all of the orbits do not converge.
It has been suggested that such
orbits could be classified using higher order averages. In the case of a bounded
observable, we show that a classical result
of G.H.~Hardy implies that if the Birkhoff averages do not converge, then neither do
the higher order averages.

If the Birkhoff averages do not converge then we may denote
by $[\alpha_k,\beta_k]$ the limit set of  the $k$-th order averages.
The sequence of intervals thus generated is nested:
$  [\alpha_{k+1},\beta_{k+1}] \subset [\alpha_k,\beta_k]$.
We can thus make a distinction among nonconvergent
Birkhoff averages; either:
\begin{enumerate}
\item[\bf B1.] {\em $\cap_{k=1}^\infty [\alpha_k,\beta_k]$ is a point
               $B_\infty$}, or,
\item[\bf B2.] {\em $\cap_{k=1}^\infty [\alpha_k,\beta_k]$ is an interval
                     $[\alpha_\infty,\beta_\infty]$.}
\end{enumerate}
We give characterizations of the types {\bf B1} and {\bf B2} in terms of how slowly
they oscillate and we give examples that exhibit both behaviours {\bf B1} and {\bf B2}
in the context of full shifts on finite symbols and ``Bowen's example''.
For finite full shifts,
we show that the set of orbits with type {\bf B2} behaviour has full topological entropy.
\end{abstract}

\section{Birkhoff averages and higher order averages}

Let $X$ be a metric space, let $f$ be a continuous map from $X$ into itself
and let $\phi:X \rightarrow \RR$ be a  continuous function (observable).  Usually
we will assume that $\phi$ is bounded.
Given $x \in X$, denote $x_i = f^i(x)$ and consider the sequence  of partial means:
\begin{equation}\label{eq:ba}
    B_n(x) \equiv \frac{1}{n} \sum_{i=0}^{n-1} \phi(x_i).
\end{equation}
The limit of this sequence, if it exists, is called the
Birkhoff average of $\phi$ along the orbit $\{x_i\}$. We can
also consider averages along orbits of continuous flows on $X$ in
the obvious way. The limit (\ref{eq:ba}) need not exist, in which case the limit set
of the sequence $\{B_n\}$ is a closed interval in $\RR$ which we
will denote by $[\alpha_0,\beta_0]$. Such orbits have been labelled {\em historical}
\cite{Tak2}.

There are well-known examples of dynamical systems with orbits
whose Birkhoff averages do not exist. Among them are
``Bowen's example"  and full shifts on finite symbols, both of which we will describe later.
In Bowen's example there is an open set of initial conditions for which the
average does not converge \cite{Tak}. For full shifts, the set of orbits for
which the Birkhoff average does not converge has full Hausdorff dimension \cite{BS}.

It has been proposed that in order to study cases where the average does
not converge, one might consider the higher order averages \cite[p. 11]{BDV}.
(See the next section for definitions.)
It is suggested that they might provide a stratification of dynamical
systems or orbits of such, indicating their statistical complexity.
However, we find as a corollary to a result by Hardy \cite{Har}
that if $\phi(x_i)$ is bounded, then the higher-order
averages are convergent if and only if the Birkhoff partial averages $\{B_n\}$ converge.
Thus, if $\phi(x_i)$ is bounded then either
\begin{enumerate}
\item[\bf A.] {\bf The Birkhoff average exists}, or
\item[\bf B.] {\bf The averages of all orders diverge.}
\end{enumerate}
In other words, a stratification by higher order means does not exist for the case
of bounded observables.

There is however a possible distinction among dynamical systems or
orbits of class {\bf B} (historical) which we now describe. In case {\bf B} we denote
by $[\alpha_k,\beta_k]$ the limit set of  $k$-th order averages. It follows easily that
the sequence of intervals thus generated is nested:
$$
  [\alpha_{k+1},\beta_{k+1}] \subset [\alpha_k,\beta_k].
$$
Therefore, either
\begin{enumerate}
\item[\bf B1.] {\em $\cap_{k=1}^\infty [\alpha_k,\beta_k]$ is a point
               $\bar{B}_\infty$}, or,
\item[\bf B2.] {\em $\cap_{k=1}^\infty [\alpha_k,\beta_k]$ is a non-trivial interval
                     $[\alpha_\infty,\beta_\infty]$.}
\end{enumerate}

We will give examples that exhibit both behaviours {\bf B1} and {\bf B2}.
Specifically, we will show that a full shift on finite symbols has orbits
of each type. In fact, for any finite full shift, we show that the set
of orbits with type {\bf B2} behaviour has full topological entropy.
We will show that all the orbits in Bowen's example are of type {\bf B1},
while in a non-hyperbolic modification of the example all orbits are of type {\bf B2}.


\section{Ces\'{a}ro and H\"{o}lder Means and Slowly Oscillating Series}

We begin with the Ces\'{a}ro means.
Suppose $\{a_n\}$ is a sequence of real numbers and denote
\begin{equation*}
 \begin{split}
   S_n &= a_1 + a_2 + \ldots + a_n \\
   S_n^{(1)} &= S_1 + S_2 + \ldots + S_n \\
    \vdots \quad &= \quad \vdots  \\
    S_n^{(k)} &= S_1^{(k-1)} + S_2^{(k-1)} + \ldots + S_n^{(k-1)} \\
   \cdots &=  \cdots + \cdots + \cdots  .
 \end{split}
\end{equation*}
We say that $\sum a_n$ is ``summable by Ces\'{a}ro mean of the $k$-th order'',
denoted $(C,k)$, if
$$
   C^{(k)}_n \equiv  \frac{S_n^{(k)}}{n^k}
$$
converges as $n \rightarrow \infty$. In \cite{Har} we find the following result.
\begin{thm}[Hardy]
If $|n a_n| < K$ for all $n \ge 1$, then the series $\sum a_n$ cannot be summable by Ces\'{a}ro mean
of any order unless it is convergent.
\end{thm}
An explanation of this Theorem is that a divergent series under the condition
$|n a_n|<K$ must oscillate very slowly, and this slow oscillation
also occurs in the Ces\'{a}ro means.

In the context of Birkhoff averages, if we let $a_1 = \phi(x_0)$ and
$$
a_n = \frac{1}{n} \phi(x_{n-1}) - \frac{1}{n} B_{n-1}, \quad \textrm{ for } \quad n>1,
$$
and define Ces\'{a}ro means as above, then we have $S_n(x) = B_n(x)$.
Note that if $\{\phi(x_i)\}$ is bounded, then $|n a_n| \le 2 \sup_i|\phi(x_i)|$ and so
Hardy's Theorem has the following consequence in the context of Birkhoff averages.
\begin{cor}\label{cor:div}
Suppose $\phi(x_i)$ is bounded and the Birkhoff average of
$\phi$ along $\{x_i\}$ diverges, then the Ces\'{a}ro means of all orders
of the Birkhoff averages also diverge.
\end{cor}

In \cite[p. 11]{BDV} higher order means by the method of H\"{o}lder was suggested.
For $\{S_n\}$ as above define
\begin{equation*}
      H^{(1)}_n = \frac{1}{n} \sum_{i=1}^{n} S_i,
\end{equation*}
and define recursively
\begin{equation*}
     H^{(k)}_n = \frac{1}{n} \sum_{i=1}^{n} H^{(k-1)}_i.
\end{equation*}
The series $\sum a_n$ is said to be {\em summable by H\"{o}lder means} or
{\em summable} (H,k) if the sequence $\{H^{(k)}_n\}$
converges as $n \rightarrow \infty$. Note that $H^{(1)}_n = C^{(1)}_n$, so first order Ces\'{a}ro and H\"{o}lder means coincide. For $k>1$ the means differ, but the methods of summation are
equivalent in the sense that $\sum a_n$ is summable $(C,k)$ if and only if it
is summable $(H,k)$. (This result can be found in \cite[\S 5.8]{HarBook} based on  \cite{Kno} and \cite{Sch}.) Thus the conclusions of Corollary~\ref{cor:div} also hold for
H\"{o}lder means of all orders.

Below we will denote $H^{(0)}_n \equiv B_n = S_n$, and, $H^{(k)}_n$, $n \ge 1$ will denote
the $k$-th order H\"{o}lder means of $B_n$.

As pointed out by Hardy, divergent sums with the bound $|n a_n| <K$ must
oscillate slowly. We give a characterization of the slowness.
Suppose that $\{B_n\}$ is divergent, with limit set $[\alpha_0,\beta_0]$.
Given any $\ep>0$ the sequence $\{B_n\}$ is infinitely often
within $\ep$ of each of the endpoints $\alpha_0$ and $\beta_0$.
For $\ep >0$, define
$$
    t_1(\ep) = \{ \min n \ge 1 : B_n(x) > \beta_0 - \ep \},
$$
and
$$
t_2(\ep) = \{ \min n > t_1 :  B_n(x) < \alpha_0 + \ep \}.
$$
If $\ep$ is sufficiently small, we may define a unique sequence $\{t_j\}$ by
\begin{equation*}
  t_{j}(\ep) = \min n > t_{j-1}(\ep) :
     \begin{cases} B_n > \beta_0 -\ep , \textrm{ if } j \textrm{ is odd }\\
                      B_n < \alpha_0 + \ep , \textrm{ if } j \textrm{ is even }.
        \end{cases}
\end{equation*}

\begin{prop}\label{expon}
Suppose $\phi(x_i)$ is bounded. If $\{B_n(x)\}$ is divergent, then for any sufficiently small
$\ep >0$ there exists $d>1$ such that
$$
     \frac{t_{j+1}(\ep)}{t_j(\ep)} \ge d, \qquad \textrm{ for all } \quad j \ge 1.
$$
\end{prop}
\noindent
{\bf Proof:}\\
We will denote
$$
\alpha_{-1} = \liminf_{i \rightarrow \infty} \phi(x_i) \qquad \textrm{ and }
       \qquad \beta_{-1} = \limsup_{i \rightarrow \infty} \phi(x_i),
$$
Let $\ep$ be small enough so
that $\beta_0 - \ep > \alpha_0 +\ep$ and assume that $n$ is large enough so that
$\alpha_{-1}-\ep < \phi(x_i) < \beta_{-1} +\ep$ for all $i\ge n$. In the following
we drop the dependence of $t_j$ on $\epsilon$.

Suppose $j$ is even so that
$$
    B_{t_j}(x) < \alpha_0 + \ep.
$$
We have then that
\begin{equation*}
 \begin{split}
    B_{t_j + 1} &< \frac{1}{t_j + 1} ( t_j B_{t_j}(x) + \beta_{-1} +\ep ) \\
          &< \frac{t_j}{t_j + 1} (\alpha_0+\ep) + \frac{\beta_{-1} +\ep}{t_j + 1},
 \end{split}
\end{equation*}
and for any $i \ge 1$,
\begin{equation}
 B_{t_j + i}< \frac{t_j}{t_j + i} (\alpha_0+\ep) + \frac{i(\beta_{-1} +\ep)}{t_j + i}.
\end{equation}
Since
$$
     B_{t_{j+1}} > \beta_0 - \ep,
$$
we must have
$$
  \frac{1}{t_{j+1}} \left( t_j (\alpha_0 + \ep) + (t_{j+1} - t_j) (\beta_{-1} +\ep) \right)
          > \beta_0 - \ep.
$$
Solving we find
$$
    \frac{t_{j+1}}{t_j} > \frac{\beta_{-1} - \alpha_0 }{\beta_{-1} - \beta_0 + 2 \ep} >1.
$$
Similarly, for $j$ odd we obtain
$$
    \frac{t_{j+1}}{t_j} > \frac{\beta_0 -\alpha_{-1}}{\alpha_0 - \alpha_{-1} + 2\ep} > 1,
$$
and thus we can set
$$
    d = \min\left\{ \frac{\beta_0 - \alpha_{-1}}{\alpha_0 - \alpha_{-1} + 2 \ep}
               , \frac{\beta_{-1}  - \alpha_0 }{\beta_{-1}  - \beta_0 + 2 \ep} \right\}.
$$
\hfill $\Box$

Thus if a Birkhoff average diverges it must oscillate with at least
exponentially increasing times. We now
distinguish between those orbits which will exhibit behaviour
\textbf{B1} and those which will exhibit behaviour \textbf{B2}.
To this end for any $\gamma \in [\alpha_0,\beta_0]$
and $\epsilon>0$ we will let $n_i = n_i(\gamma,\epsilon)$  be the
subsequence of positive integers such that
$B_{n_i} \in (\gamma-\epsilon,\gamma+\epsilon)$. Our main result will
be the following.
\begin{thm}\label{main}
Suppose that $B_n(x)$ has limit set $[\alpha_0,\beta_0]$.
\begin{enumerate}
\item[1.]
If there exists $D>1$ such that for all $\gamma \in [\alpha_0,\beta_0]$ and
$\epsilon>0$
\begin{equation}\label{n_i-bdd}
\frac{n_{i+1}(\gamma,\epsilon)}{n_i(\gamma,\epsilon)} \leq D
    \text{ for all }i\in\mathbb{N}
\end{equation}
then $B_n(x)$ is \textbf{B1}.
\item[2.]
If for every $\gamma \in [\alpha_0,\beta_0]$
we can find $\epsilon>0$ such that
$$
\limsup_{i\rightarrow\infty}\frac{n_{i+1}(\gamma,\epsilon)}{n_i(\gamma,\epsilon)}=\infty
$$
then $B_n(x)$ is \textbf{B2}
\end{enumerate}
\end{thm}


\section{Proof of Theorem~\ref{main} and Other Results}

\begin{prop}\label{prop:contraction1}
Suppose that $\phi(x_i)$ is bounded, $B_n$ diverges and there exists $1<D<\infty$ such that
\begin{equation}\label{limsup}
    \limsup_{j \rightarrow \infty} \frac{ t_{j+1}(\epsilon)}{t_j(\epsilon)} \le D,
\end{equation}
for all $\ep >0$ sufficiently small, then
\begin{equation}\label{eqn:contraction1}
   \beta_0 - \alpha_0 \le \frac{D-1}{D+1} \left(\beta_{-1} - \alpha_{-1}\right).
\end{equation}
\end{prop}
\noindent
{\bf Proof:}\\
Denote $\alpha_{-1}$ and $\beta_{-1}$ as above and assume that $n$ is large
enough so that $\alpha_{-1}-\ep < \phi(x_i) < \beta_{-1} +\ep$ for all $i \ge n$.
Suppose $j$ is even so that $B_{t_j} < \alpha_0 + \ep$. Given
any $\delta >0$ we may assume that $j$ is large enough so that $t_{j+1}/t_j < D + \delta$.
Then we have
\begin{equation}
 \begin{split}
  B_{t_{j+1}} - B_{t_j} &\le \frac{1}{t_{j+1}}
                 \left( t_j B_{t_j} + (t_{j+1} - t_j)(\beta_{-1} + \ep) \right) - B_{t_j} \\
              &\le \frac{t_{j+1} - t_j}{t_{j+1}} \left( \beta_{-1} + \ep - B_{t_j} \right) \\
              & \le \left( 1 -  \frac{1}{D+\delta} \right)
                       \left( \beta_{-1}  - \alpha_0 \right).
 \end{split}
\end{equation}
Thus we conclude that
$$
\beta_0 - \alpha_0 \le \frac{D-1}{D} (\beta_{-1} - \alpha_0).
$$
Similarly, for $j$ odd we calculate:
$$
   \beta_0 - \alpha_0 \le \frac{D-1}{D} (\beta_0 - \alpha_{-1}).
$$
Combining these inequalities we obtain (\ref{eqn:contraction1}).
\hfill $\Box$

Now consider higher order means. For $\{B_n\}$ we defined a sequence of
times $\{t_j\}$ above. If $j$ is even, recall that $B_{t_j} < \alpha_0 +\ep$
and $B_{t_{j+1}} > \beta_0 - \ep$.  Define $t^{(1)}_j$ to be the index $n$,
$t_j \le n < t_{j+1}$, at which $H^{(1)}_n$ achieves its minimum. Also for
$j$ odd define $t^{(1)}_j$ to be the index $n$,
$t_j \le n < t_{j+1}$, at which $H^{(1)}_n$ achieves its maximum.
Note that $H^{(1)}_n$ can achieve a local minimum or maximum only
when its values and those of $B_n$ cross.

For $k>1$ we may define inductively times $\{t^{(k)}_j\}$ in a similar way as $\{t^{(1)}_j\}$.
In particular we have:
\begin{equation}\label{eqn:times}
\cdots < t^{(k-1)}_{j-1}  \le   t^{(k)}_{j-1} <  t^{(k-1)}_j \le  t^{(k)}_j
                <  t^{(k-1)}_{j+1}  \le  t^{(k)}_{j+1}  < \cdots.
\end{equation}

The proof of the following proposition is similar to that of Proposition~\ref{prop:contraction1}.
\begin{prop}\label{prop:contraction2}
Suppose that $\phi(x_i)$ is bounded, $B_n$ diverges and there exists $1<D<\infty$ such that
\begin{equation}\label{limsup2}
    \limsup_{j \rightarrow \infty} \frac{ t^{(k)}_{j+1}}{t^{(k)}_j} \le D,
\end{equation}
for some $k$ and for all $\ep >0$ sufficiently small, then
\begin{equation}\label{eqn:contraction2}
   \beta_k - \alpha_k \le \frac{D-1}{D+1} \left(\beta_{k-1} - \alpha_{k-1}\right).
\end{equation}
\end{prop}
Note however, that (\ref{limsup}) does not imply (\ref{limsup2}), but the nesting
(\ref{eqn:times}) along with (\ref{limsup}) imply that the limit supremum
in (\ref{limsup2}) is less than $D^{k+1}$.

Next we consider type {\bf B2} sequences.
A specific case of type {\bf B2} behavior is captured in the following.
\begin{prop}\label{prop:nondecreasing}
If $\phi(x_i)$ is bounded, $\{B_n\}$ diverges and
$[\alpha_0, \beta_0] = [\alpha_{-1},\beta_{-1}]$
then $[\alpha_k,\beta_k] = [\alpha_{-1},\beta_{-1}]$ for all $k$,
so the Birkhoff averages are of type {\bf B2}.
\end{prop}
\noindent
{\bf Proof:}\\
Since the limit set of $\{H^{(k)}\}$ is $[\alpha_k,\beta_k]$
it follows that for a given $k$ and any $\ep>0$ there exists
$N(k,\ep)$ such that $\alpha_k-\ep < H^{(k)}_n < \beta_k+\ep$ for all $n>N(k,\ep)$.

Under the assumptions of the proposition, given any $\ep>0$, the sequence
$B_n$ is greater than $\beta_{-1} - \ep$ infinitely often. Further, we may assume
that $n$ is sufficiently large so that $\phi(x_n) < \beta_{-1} + \ep$.  Now given any
$\delta > \ep$ suppose that:
$$
   B_n < \beta_{-1} - \delta  \quad \textrm{while} \quad B_{n+i} > \beta_{-1} - \ep,
$$
for some positive integer $i$. By the assumptions this must happen infinitely often.
Since $\phi(x_n) < \beta_{-1} + \ep$ it follows that
$$
B_ {n+i} < \frac{1}{n+i} \left( n (\beta_{-1} - \delta) + i (\beta_{-1} + \ep) \right).
$$
Since  $B_{n+i} > \beta_{-1} - \ep$ we find from the preceding equation that
$$
  i >  \frac{n(\delta - \ep)}{2\ep}.
$$
If $n$ is the last integer before $B_j$ reaches a maximum but $B_n < \beta_{-1} - \delta $, then
$B_j > \beta_{-1} -\delta$ for at least $i = n(\delta -\ep)/2\ep$ consecutive steps. This gives us that
\begin{equation}
\begin{split}
H^{(1)}_{n+i} &=  \frac{1}{n+i} \left( n H^{(1)}_n + \sum_{j = n+1}^{n+i} B_j \right) \\
   &  >  \frac{1}{n+i} \left( n (\alpha_{-1} -\ep) + i ( \beta_{-1}-\delta) \right) \\
   & > \frac{2\ep}{\delta+\ep} (\alpha_{-1} -\ep)
             + \frac{\delta-\ep}{\delta + \ep} (\beta_{-1} - \delta).
\end{split}
\end{equation}
Since $\ep$ is taken arbitrarily small (by taking $n$ large) $H^{(1)}_{n+i}$ is arbitrarily close to $\beta_{-1} - \delta$.
Since $\delta$ is taken arbitrarily small, we conclude that $\beta_1 = \beta_{-1}$.

Assuming that $\overline{\lim}_{n \rightarrow \infty} H^{(k)}_n = \beta_{-1}$, by a similar calculation as above we can show that
$\overline{\lim}_{n \rightarrow \infty} H^{(k+1)}_n = \beta_{-1}$, which proves the result.
\hfill $\Box$

Now we prove part 1 of Theorem~\ref{main}.
First note that if (\ref{n_i-bdd}) holds for both $\gamma = \alpha_0$ and $\gamma = \beta_0$,
then all $j$ we have $\frac{t_{j+1}(\epsilon)}{t_j(\epsilon)}\leq D$.

\begin{lem}\label{lem:b2}
 If $\phi(x_i)$ is bounded, $\{B_n\}$ diverges and
the Birkhoff averages are of type {\bf B2}, then for any $k \ge 1$
\begin{equation}\label{tkratio}
     \limsup_{j \rightarrow \infty} \frac{t^{(k)}_{j+1}}{t^{(k)}_j} = + \infty.
\end{equation}
\end{lem}
\noindent
{\bf Proof:}
By the assumptions there exists $K$ such that
$$
\alpha_k > \alpha_\infty - \frac{\ep}{2}   \quad \textrm{and}  \quad  \beta_k < \beta_\infty + \frac{\ep}{2},
$$
for all $k \ge K$. Note that $H^{(k)}_n > \beta_k - \ep/2$ for infinitely many $n$.
Given $\delta > \ep$ suppose that
$$
   H^{(k)}_n < \beta_k - \delta  \quad \textrm{and}
                 \quad H^{(k)}_{n+i} > \beta_k - \frac{\ep}{2}.
$$
Since $H^{(k-1)}_n < \beta_\infty + \ep$ for $k > K $ and $n$ sufficiently large, we see that
\begin{equation*}
H^{(k)}_{n+i} < \frac{1}{n+i} \left( n (\beta_k - \delta) + i (\beta_\infty + \ep) \right).
\end{equation*}
In order to have $H^{(k)}_{n+i} > \beta_k - \frac{\ep}{2}$
we calculate that $i$ must satisfy:
$$
      i > \frac{n(\delta - \ep)}{2\ep}.
$$
Let $t^{(k)}_j$ be the minimum point before $n$ and $t^{(k)}_{j+1}$ be the maximum following $n$.
We then have that
$$
       \frac{t^{(k)}_{j+1}}{t^{(k)}_j} > \frac{\delta-\ep}{2\ep}.
$$
Since $\ep$ can be taken arbitrarily small the result follows for all $k >K$.
The interlacing (\ref{eqn:times}) of $\{t^{(k-1)}\}$ with $\{t^{(k)}\}$
implies the result for $k \le K$.
\hfill $\Box$

This gives us a sufficient condition for {\bf B1} behaviour.
\begin{cor}\label{cor:b1}
Suppose that $\phi(x_i)$ is bounded, $B_n$ diverges and there exists $1<D<\infty$ such that
\begin{equation}\label{limsup4}
    \limsup_{j \rightarrow \infty} \frac{ t_{j+1}}{t_j} \le D,
\end{equation}
for all $\ep >0$ sufficiently small, then the Birkhoff averages
are of type {\bf B1}.
\end{cor}
\noindent
{\bf Proof:}\\
Given the nested definition of the times in (\ref{eqn:times}), we can conclude
that
\begin{equation}\label{t1ratio}
    \limsup_{j \rightarrow \infty}  \frac{t^{(1)}_{j+1}}{t^{(1)}_j} < D^2.
\end{equation}
The result then follows from Lemma~\ref{lem:b2}.
\hfill $\Box$

We can deduce part 1 of Theorem~\ref{main} from Corollary~\ref{cor:b1}.

Next, let us consider part 2 of the Theorem.
Recall our notation that $[\alpha_0,\beta_0]$ is the limit set for the
sequence of Birkhoff averages $B_n$.
\begin{lem}
For any $\gamma\in [\alpha_0,\beta_0]$
if there exists an $\epsilon>0$ such that
\begin{equation}\label{limsup3}
\limsup_{i\rightarrow\infty}\frac{n_{i+1}(\epsilon,\gamma)}{n_i(\epsilon,\gamma)}=\infty,
\end{equation}
then $\cap_{k=1}^{\infty}[\alpha_k,\beta_k]$ contains a point outside of $\left(\gamma-\frac{\epsilon}{2},\gamma+\frac{\epsilon}{2}\right)$. In  particular
$\cap_{k=1}^{\infty}[\alpha_k,\beta_k]$ cannot be the single point $\gamma$.
\end{lem}
\begin{proof}
Fix $\gamma \in [\alpha_0,\beta_0]$.
We define a sequence $\{t_i\}_{i\in N}$ as follows
$$t_1 \epsilon -(\gamma - \alpha_0 +\epsilon) \geq \frac{(1+t_1)3\epsilon}{4}$$
and for $n> 1$
$$
 t_n \, \epsilon \, \frac{2^{n-1}+1}{2^n}-\left(1+\sum_{i=1}^{n-1}t_i\right)\alpha_0
 \geq \frac{2^n+1}{2^{n+1}}\left(1+\sum_{i=1}^n t_i\right) n \, \epsilon.
$$
Fix $k\in \NN$. Since \ref{limsup3} holds we can find infinitely many $n_i$ such that $\frac{n_{i+1}}{n_i}\geq 1+\sum_{j=1}^k t_j$. We now show that for any $1\leq n\leq k$,
the interval $[\alpha_n,\beta_n]$ must contain a point outside  $\left(\gamma-\frac{\epsilon}{2},\gamma+\frac{\epsilon}{2}\right)$. For convenience
we will assume that for each integer $i$ and for each $n$ between $n_i$ and $n_{i+1}$,
we have $B_n\geq\alpha+\epsilon$ (To handle the other case we switch inequalities and
replace $\alpha_0$ by $\beta_0$).
We proceed by induction.
To start with consider $H_j^{(1)}$ for $ (n_i)(1+t_1)\leq j< n_{i+1}$ we know that
$$
jH_j^{(1)}-j\gamma
   \geq (t_1)(n_i)\epsilon+(j-n_i(1+t_1))\epsilon-n_i(\gamma - \alpha_0 +\epsilon)
     >  \frac{n_i(1+t_1)3\epsilon}{4}+(j-n_i(1+t_1))\epsilon
$$
from which it follows that $H_j^{(1)} \geq \gamma +\frac{3\epsilon}{4}$.
Assume that for $ (n_i)\left(1+\sum_{l=1}^{n-1}t_l\right) \leq j<n_{i+1}$
we have $h_{j}^{(n-1)}-\gamma \geq \frac{2^{n-1}+1}{2^n} \, \epsilon$ thus
for $z \geq (n_i) \left( 1+\sum_{l=1}^{n}t_l \right)$ we know that
\begin{eqnarray*}
zH_z^{(n)}-z\gamma
   &\geq& \left(z - n_i \left( 1+\sum_{l=1}^{n-1}t_l\right ) \right)
    \frac{2^{n-1}+1}{2^n} \, \epsilon - n_i (\gamma - \alpha_0 + \epsilon)
       \left( 1+\sum_{l=1}^{n-1}t_l \right)   \\
&=& n_i t_n \epsilon \frac{2^{n-1}+1}{2^n} - \alpha_0 n_i
    \left( 1+\sum_{l=1}^{n-1}t_l \right)
  + \left(z-n_i \left(1+\sum_{l=1}^n t_l \right) \right)
       \frac{2^{n-1}+1}{2^n} \,\epsilon\\
&\geq& \frac{2^n+1}{2^{n+1}} \left( 1+\sum_{l=1}^n t_l \right)
     n_i \epsilon + \left(z-n_i\left(1+\sum_{l=1}^n t_l\right)\right)
         \frac{2^{n-1}+1}{2^n} \, \epsilon
\end{eqnarray*}
and it follows that $H_z^{(n)}\geq \gamma + \frac{2^n+1}{2^{n+1}}$.
The result now follows by induction.

\end{proof}
Part 2 of Theorem~\ref{main} immediately follows.


\section{Full Shift on Finite Symbols}

\subsection{Some Examples}

{\bf Example 1.} {\em A sequence whose orbit is type {\bf B1}}.\\
Consider the sequence starting with $0$, followed by 2 ones,
followed by $2^2$ zeros, followed by $2^3$ ones, etc.., i.e.
$$
\{x_n\} = \{0, 1, 1,  0, 0, 0, 0, \underset{2^3}{\underbrace{1, \ldots, 1}},
           \underset{2^4}{\underbrace{0, \ldots, 0}}, 1, \ldots \}.
$$
For this example, one may calculate directly that the limit set of $B_n$ is
$[1/3,2/3]$.
Let $t_j$ be as above, then
$$
   \lim_{j \rightarrow \infty}  \frac{t_{j+1}}{t_j} = 2 = d = D.
$$
Thus, by Corollary~\ref{cor:b1} the orbit of this sequence is of type {\bf B1}.
We also observe that the estimate (\ref{eqn:contraction1}) is optimal.

\smallskip
\noindent
{\bf Example 2.} {\em A sequence whose orbit is of type {\bf B2}.}\\
Consider the sequence beginning with 1 zero, followed by 2 ones,
followed by 9 zeros, followed by 48 ones, etc., so that
the $i$th group of constant $x_n$ is $i$ times as long as
all the preceding groups put together. That is, the length $\ell_i$ of the
$i$-th constant group is given recursively as:
$$
  \ell_1 = 1,  \qquad  \ell_i = i \sum_{k=1}^{i-1} \ell_k, \textrm{ for } i \ge 2.
$$
It is easily seen for this example that the limit set of $B_n$ is
$[0,1]$. By Proposition~\ref{prop:nondecreasing} the sequence is of type {\bf B2}.

\smallskip
\noindent
{\bf Example 3.} {\em An $L^1$ counterexample.}
We can construct an unbounded $L^1$ function and a
sequence for which the first Birkhoff averages do not
converge but the second averages converge.
Consider the space $\{-1,1\}^\NN$ and let $\sigma$ be the usual left shift. We define $f_1,f_2,f:\Sigma\rightarrow\RR$ as follows
$$f_1(i)=\inf\{n:i_{n+1} \neq i_1\},    $$
$$f_2(i)=i_1,$$
$$f=f_1f_2. $$
$f$ is clearly $L^1$ with respect to any $(p,1-p)$ Bernoulli measure but is not
continuous since it blows up at $(1,1,1,\ldots)$ and $(-1,-1,\ldots)$.
\begin{prop}
For $f$  and $j=(1,-1,-1,1,1,1,-1,-1,-1,-1,\ldots)$ the Birkhoff averages
do not converge. However the average of the Birkhoff averages does converge.
\end{prop}
\noindent
{\bf Proof:}
Let $k=n(n+1)/2$. If $n$ is odd then
$$
\sum_{m=0}^{k-1}f(\sigma^m j)= 1-(1+2)+(1+2+3)-\ldots+(1+2+3+\ldots+n)=1+3+5+\ldots+n.
$$
If $n$ is even then
$$
\sum_{m=0}^{k-1}f(\sigma^m j)= 1-(1+2)+(1+2+3)-\ldots-(1+2+3+\ldots+n)=-2-4-6-\ldots-n.
$$
from which we can clearly deduce that the Birkhoff averages from $f$ at $j$ do
not converge (note this does not use exponential times of oscillation but instead the fact
that $f$ is unbounded). In fact the averages oscillate between $-\frac{1}{2}$ when
$n$ is even and $\frac{1}{2}$ when $n$ is odd.

Now we need to show that the average of the Birkhoff averages does converge.
 It is clear from the first part that the averages oscillate between $-\frac{1}{2}$ and
$\frac{1}{2}$ subexponentially. By Proposition~\ref{expon} the average of this sequence
must converge.

\hfill $\Box$



\subsection{Topological entropy}

In a finite shift the set for which the
averages do not converge is known to have full topological entropy or
alternatively Hausdorff dimension (excluding the case when the
function is cohomologous to a constant). A definition of topological
entropy for non-compact sets was introduced in \cite{Bowen}. For further
discussions see \cite{TV} and \cite{BS}.

Let $\Sigma$ be the full shift on $m$ symbols and
$f:\Sigma\rightarrow\mathbb{R}$ a continuous function which
is not cohomologous to a constant. We will show that the set
of points for which the Birkhoff averages are in class {\bf B2}
also has full topological entropy.

\begin{prop}
The set of points in class {\bf B2} for $f$ has topological entropy equal to $\log m$.
\end{prop}
\begin{proof}
Let $\nu$ be the evenly weighted Bernoulli measure and $\alpha=\int f\text{d}\nu$.
Fix $\alpha_1\neq\alpha_2\in\mathbb{R}$ such that we can find two ergodic shift invariant probability measures $\mu_1,\mu_2$ such that $\int f\text{d}\mu_1=\alpha_1$ and $\int f\text{d}\mu_2=\alpha_2$. Choose $0<\epsilon<\frac{\alpha_2-\alpha_1}{4}.$
Let $X_N$ consist of sequences $\omega$ such that for all $n\geq N$
$$\left|\sum_{i=0}^{n-1}f(\sigma^i\omega)-n\alpha_1\right|\leq n\epsilon$$
and
$$\left|\mu_1([i_1,\ldots,i_n])-nh(\mu_1)\right|\leq n\epsilon.$$
Similarly, let $Y_n$ consist of sequences $\omega$ such that
$$\left|\sum_{i=0}^{n-1}f(\sigma^i\omega)-n\alpha_2\right|\leq n\epsilon$$
and
$$\left|\mu_2([i_1,\ldots,i_n])-nh(\mu_2)\right|\leq n\epsilon.$$
Note that by the Birkhoff Ergodic Theorem and the Shannon-McMillan-Brieman
Theorem~\cite[p. 93]{Wa}
$$
\lim_{N\rightarrow\infty}\mu_1(X_N)=1\text{ and }\lim_{N\rightarrow\infty}\mu_2(Y_N)=1.
$$
We now construct a new subset $Z_N$.
We let $n_1=N$ and $n_i=i\sum_{j=1}^{i-1}n_j$ for $i\geq 2$. We also let
$k_i=\sum_{j\text{ odd}}^{2i-1}n_i$, $l_i=\sum_{i\text{ even }}^{2i}n_i$
be the sum of the odd and even $n_i$s and for convenience let $k_0=l_0=0$.
We define $Z_n$ by the condition a sequence $\tau\in Z_n$ if and only if for each $i\geq 0$
\begin{enumerate}
\item
 $[\tau_{k_i+l_i+1},\ldots,\tau_{k_i+l_i+n_{2i+1}}]\cap X_N\neq\emptyset$
\item
$[\tau_{k_{i+1}+l_i+1},\ldots,\tau_{k_{i+1}+l_{i}+n_{2i+2}}]\cap Y_N\neq\emptyset$.
\end{enumerate}
Using the continuity of $f$ and the definition of $Z_n$ we can see that
for any $\tau\in Z_n$ the limit set of the Birkhoff averages contains the
interval $[\alpha_1+\epsilon,\alpha_2-\epsilon]$ and the condition for
Proposition~\ref{prop:nondecreasing}  are clearly satisfied.
We now fix $N$ large enough so that $\mu_1(X_N),\mu_2(Y_N)\geq\frac{1}{2}.$ Thus we can estimate the growth of the number of $n$th level cylinders in $Z_N$. It is clear that for $n\geq N$
\begin{eqnarray*}
\#\{[[\tau_1,\ldots,\tau_n]:[\tau_1,\ldots,\tau_{n}]\cap X_N\neq\emptyset\}&\geq&\frac{1}{2}e^{n(h(\mu_1)-\epsilon)}   \\
\#\{[[\tau_1,\ldots,\tau_n]:[\tau_1,\ldots,\tau_{n}]\cap Y_N\neq\emptyset\}&\geq&\frac{1}{2}e^{n(h(\mu_2)-\epsilon)}.
\end{eqnarray*}
Combining this with the definition of $Z_n$ we can see that if we let $M_n$ be the number of $n$th level cylinders containing elements of $Z_n$ then
$$
\liminf_{n\rightarrow\infty} \frac{\log{M_n}}{n}
    \geq  \min \{h(\mu_1)-\epsilon,h(\mu_2)-\epsilon\}.
$$

Thus we can define a measure $\nu$ by
$$
\nu([\tau_1,\ldots,\tau_n])
    =   \left\{\begin{array}{lll}
           \frac{1}{M_n} & \text{ if } &[\tau_1,\ldots,\tau_n]\cap Z_N\neq\emptyset\\
           0 & \text{ if } & [\tau_1,\ldots,\tau_n]\cap Z_N=\emptyset
         \end{array}\right.
$$
This measure will satisfy $\nu(Z_N)=1$ and for any $\tau\in Z_N$ and $n$ sufficiently large
$$
\nu([\tau_1,\ldots,\tau_n])\leq e^{-n\min\{h(\mu_1)-\epsilon,h(\mu_2)-\epsilon\}}.
$$
It follows by the entropy distribution principle (Theorem 3.6 in \cite{TV}) that
$$
H(Z_n) \geq \min\{h(\mu_1)-\epsilon,h(\mu_2)-\epsilon\}.
$$
To complete the proof we note that using results in \cite{TV} for any $\delta>0$,
by choosing $\alpha_1$ and $\alpha_2$ sufficiently close to $\alpha$  we can find
measures $\mu_1$ and $\mu_2$ where $\log m-h(\mu_1)<\delta$, $\log m-h(\mu_2)<\delta$
and $\int f\text{d}\mu_1=\alpha_1\neq\int f\text{d}\mu_2=\alpha_2$.
\end{proof}


\section{Bowen's example and a modification}

\subsection{Bowen's example is of type {\bf B1}.}

It is well known that in Bowen's example there is an open set of initial
conditions whose orbits are historical, i.e.~of type {\bf B}. In this
section we will recall Bowen's example and
show that its orbits are in fact of type {\bf B1}.

Let $F^t$ be a flow possessing two hyperbolic equilibria
$\mathbf{p}_1$ and $\mathbf{p}_2$ and a heteroclinic cycle as shown
in Figure~1. Let the flow be symmetric under the transformations:
$$
(u,t) \mapsto (-u,-t) \qquad \text{and} \qquad (v,t) \mapsto (-v,-t).
$$
(This symmetry is not essential, but we assume it for clarity.)
In this section we will use the Birkhoff average for the flow, which is
defined as
$$
   B(x) = \lim_{t\rightarrow +\infty}
            \frac{1}{t} \int_0^t \phi(F^t(x)),
$$
if it exists.

Suppose that $\phi(p_1) \neq \phi(p_2)$ and that the
linearized flow at each of the two hyperbolic
equilibria has eigenvalues
$$
   -\la < 0 < \mu,
$$
and suppose
$$
  \rho = \frac{\la}{\mu} > 1.
$$
so that all orbits in an interior neighbourhood of the heteroclinic
cycle limit onto the cycle.
Let  $U_1$ and $U_2$ be small neighbourhoods of these points.
It follows that the flow on these neighbourhood is $C^1$ linearizable.
In the linearized coordinate the flow has the form:
\begin{equation}
 \begin{split}
  \dot{x} &= \mu x \\
  \dot{y} &= -\lambda y
 \end{split}
\end{equation}
In fact we may choose the linearizing transformation $\Psi_i:(u,v) \mapsto (x,y)$ so that
$\Psi_i: U_i \mapsto [0,1]\times [0,1]$ \cite{Yo}.
As interior orbits circulate toward the heteroclinic cycle, they spend
more and more time near the equilibria and the average moves first
toward $\phi(p_1)$ then toward $\phi(p_2)$. To be precise, the time of transition
through $U_1$ or $U_2$ is
$$
T_j = \frac{1}{\mu} \ln \frac{1}{x_j}
$$
where $x_j$ is the $x$ coordinate at which the orbit enters one of the neighbourhoods
for the $j$-th time. The orbit will then leave the neighbourhood with $y$-coordinate
$$
y_j = x_j^\rho.
$$

\begin{figure}[ht]

\centerline{\hbox{\epsfig{figure=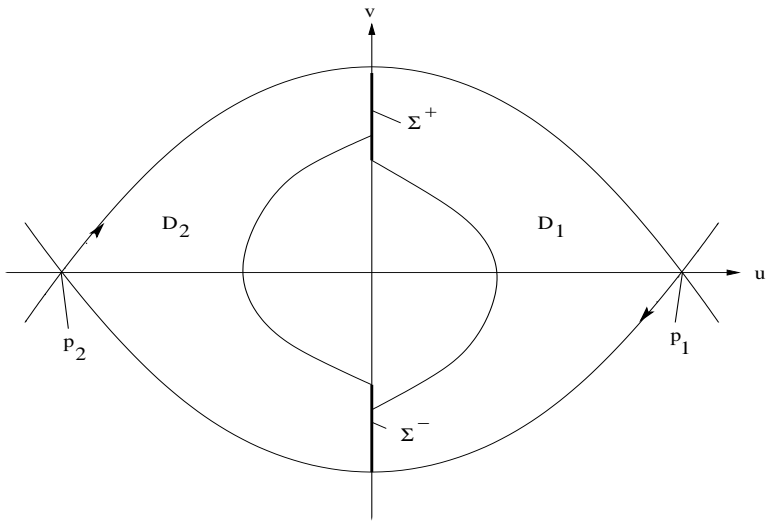,height=2.7in,width=4in}}}


\caption{A symmetric heteroclinic cycle. Cross sections $\Sigma^-$
and $\Sigma^+$ sweep out regions $D_1$ and $D_2$. They also sweep
out foliations of these regions which are used to define linearized
coordinates.}

\end{figure}

Now supposing that we begin with $x = x_0$ entering one of the neighbourhoods,
then the transition times $T_j$ are given by
$$
    T_j = C \rho^{j}
$$
where
$$
    C = \frac{1}{\mu} \ln \frac{1}{x_0}.
$$
These results also follow from Gaunerdorfer \cite{Ga} or Takens \cite{Tak}.
Thus, in Bowen's example
$$
   \lim_{j \rightarrow \infty} \frac{t_{j+1}}{t_j} = \rho,
$$
and thus it follows from Corollary~\ref{cor:b1} that the behavior is
of type {\bf B1}.
Another example involving non-hyperbolic fixed points
with the same exponential growth
as in Bowen's example was given in \cite{Yo}.

\subsection{An modification of Bowen's example with {\bf B2} behavior.}
Suppose that the equilibria, rather that being hyperbolic, have the following
form locally:
\begin{equation}
 \begin{split}
  \dot{x} &= x^3 \\
  \dot{y} &= -\lambda y,
 \end{split}
\end{equation}
on symmetric neighbourhoods in terms of charts that are given by $0 < x < d$, $0 < y < d$.
Suppose that an orbit enters one of the two neighbourhoods for the $j$-th
time at $x = x_j$ and $y = d$. The solution within the neighbourhood is
given by:
$$
     x = \frac{1}{\sqrt{x_j^{-2} - 2t}}, \qquad y = d e^{-t}.
$$
From this we calculate that the transition time $T_j$ (when $x(t) = d$) for
this crossing of the neighbourhood is
$$
     T_j = \frac{1}{2} \left( \frac{1}{x_j^2} - \frac{1}{d^2} \right).
$$
It follows that the $y$-coordinate at which the orbit leaves the neighbourhood
is:
$$
       y_j = d \exp\left(- \frac{1}{2} \left( \frac{1}{x_j^2}
                          - \frac{1}{d^2} \right) \right).
$$
The orbit will then enter the other neighbourhood at $x$-coordinate, $x_{j+1}$, given
by
$$
  x_{j+1} \approx \alpha d
          \exp\left(- \frac{1}{2} \left( \frac{1}{x_j^2} - \frac{1}{d^2} \right)\right),
$$
where $\alpha$ depends on the global flow. The corresponding transition time $T_{j+1}$
is
\begin{equation*}
 \begin{split}
     T_{j+1} &= \frac{1}{2} \left( \frac{1}{x_{j+1}^2} - \frac{1}{d^2} \right) \\
             & \approx
              \frac{ \exp\left(- \frac{1}{2} \left( \frac{1}{x_j^2}
                            - \frac{1}{d^2} \right) \right)}
                   {2 \alpha^2 d^2}
                - \frac{1}{2d^2} \\
             &  \approx  \frac{e^{T_j}}{2 \alpha^2 d^2}
 \end{split}
\end{equation*}

Thus for this example  transition times $T_j$ satisfy:
$$
   T_{j+1} \ge K e^{T_j},
$$
where $K>0$ is a constant that depends on the global flow.
Since this growth of times is even larger than in Example 2 of \S3, we
conclude that for all $k$ we have
that $[\alpha_k,\beta_k]$ is equal to the closed interval bounded by
$\phi(p_1)$ and $\phi(p_2)$. The example is thus of type {\bf B2}.

\section*{Acknowledgement}
The authors thank the Maths Institute at Warwick University
where T.J.~and V.N.~were postdoctoral fellows and T.Y.~was
a visitor when this work was conceived. T.Y.~was supported by
an Ohio University Faculty Fellowship.

\end{document}